\documentclass[leqno,draft]{article}

\begin{document}

\title{A few remarks about linear operators and disconnected open sets
in the plane}

\author{Stephen Semmes \\
        Rice University}

\date{}

\maketitle

        Let $E$ be a nonempty compact set with empty interior in the
complex plane, and let $U \subseteq {\bf C}$ be a bounded open set
such that $U \cap E = \emptyset$ and $\partial U = E$.  For instance,
$U$ might be the union of all of the bounded components of ${\bf C}
\backslash E$.  We are especially interested here in the situation
where $E$ is connected and $U$ has infinitely many connected
components, and moreover where each neighborhood of each element of
$E$ contains infinitely many components of $U$.  Of course, $U$ can
have only countably many components, since the plane is separable.
Sierpinski gaskets and carpets are basic examples of such fractal sets
$E$.

        One can choose different orientations for the different
components of $U$.  This can be represented by a locally constant
function on $U$ with values $\pm 1$, where $+1$ corresponds to the
standard orientation on ${\bf C}$.  Alternatively, a locally constant
function on ${\bf C} \backslash E$ with values $0$ and $\pm 1$ can be
used to specify which complementary components are included in $U$ as
well as their orientations.  Bergman and Hardy spaces on $U$ can then
be defined using holomorphic or conjugate-holomorphic functions on the
components of $U$, depending on the choice of orientation of the
component.

        The topological activity in $E$ indicated by its complementary
components can also be described in terms of homotopy classes of
continuous mappings from $E$ into ${\bf C} \backslash \{0\}$.  This is
reflected in Fredholm indices of associated Toeplitz operators too.
It is natural to allow different orientations on different components
of $U$, in order to get different combinations of indices, which may
involve many variations.  At the same time, the standard orientation
on all of $U$ has some special features.

        Perhaps the first point is that there are a lot of holomorphic
functions in the usual sense on neighborhoods of $\overline{U}$, which
are in particular very regular functions on $\overline{U}$ that are
holomorphic on $U$.  Conversely, removable singularity results imply
that a sufficiently well-behaved function on $\overline{U}$ which is
holomorphic on $U$ is holomorphic on the interior of $\overline{U}$.
For instance, this holds for continuously-differentiable functions on
$\overline{U}$ because $E$ has empty interior, and for Lipschitz
functions of order $1$ on $\overline{U}$ when $E$ has Lebesgue measure
$0$.  Note that the differential of a continuously-differentiable
function $f$ on $\overline{U}$ which is holomorphic on $U_1 \subseteq
U$ and conjugate-holomorphic on $U_2 \subseteq U$ vanishes on
$\overline{U_1} \cap \overline{U_2}$.  Since $E$ is connected,
$\overline{U_1} \cap \overline{U_2} \ne \emptyset$ when $U_1, U_2 \ne
\emptyset$ and $U_1 \cup U_2 = U$.

        Continuous functions that are holomorphic on both sides of a
nice curve are holomorphic across the curve.  Holomorphic functions
that are real-valued on a curve yield continuous functions that are
holomorphic on one side of the curve and conjugate-holomorphic on the
other side, by taking the complex conjugate on one side.  This is
related to the reflection principle.

        If $E$ is sufficiently big, then there are a lot of continuous
functions on ${\bf C}$ that are holomorphic on ${\bf C} \backslash E$,
given by Cauchy integrals of suitable measures on $E$.  It is not as
easy to work with the differential operator equal to
$\overline{\partial}$ on the components of $U$ with the standard
orientation and to $\partial$ on the components of $U$ with the
opposite orientation when the orientations are variable.

        Suppose that $f$ is a continuously-differentiable
complex-valued function on $\overline{U}$.  If $E$ is not
asymptotically flat at a point $p \in E$, then the differential of $f$
at $p$ is uniquely determined by the restriction of $f$ to $E$.
Complex-linearity of the differential of $f$ at $p$ is then a
significant condition, and not simply a question of choosing an
extension of $f$ on $E$ so that the differential is complex-linear, as
in the case of a smooth curve.  If $\overline{\partial} f = 0$ on $E$,
then there is more regularity for some commutators and Hankel-type
operators corresponding to $f$ and the standard orientation on all of
$U$.  There are also a lot of functions that satisfy this condition,
which is equivalent to saying that $\overline{\partial} f(p) \to 0$
when $p \in U$ approaches the boundary.  The analogous condition for
variable orientations would ask that $\overline{\partial} f(p) \to 0$
or $\partial f(p) \to 0$ as $p \in U$ approaches the boundary through
parts of $U$ with the standard or opposite orientation.  This implies
that the differential of $f$ is $0$ at any $p \in E$ in the boundary
of both parts of $U$, as before.

        Even with variable orientations on the components of $U$, one
can make use of Cauchy, conjugate-Cauchy, and other kernels on the
individual components.  Some regularity of functions on or around $E$
can still be very helpful, but not in quite the same way as when the
orientations are constant.

        Note that there is normally no uniform orientation like the
one inherited from the complex plane for fractal sets with topological
dimension $1$ bounding quasi-fractal sets with topological dimension
$2$ in ${\bf R}^n$ when $n > 2$.


\begin{thebibliography}{25}


\bibitem {a} W.~Arveson, {\it A Short Course on Spectral Theory},
Springer-Verlag, 2002.

\bibitem {b-d-f-1} L.~Brown, R.~Douglas, and P.~Fillmore, {\it
Extensions of $C^*$-algebras, operators with compact self-commutator,
and $K$-homology}, Bulletin of the American Mathematical Society {\bf
79} (1973), 973--978.

\bibitem {b-d-f-2} L.~Brown, R.~Douglas, and P.~Fillmore, {\it Unitary
equivalence modulo the compact operators and extensions of
$C^*$-algebras}, in {\it Proceedings of a Conference on Operator
Theory}, 58--128, Lecture Notes in Mathematics {\bf 345},
Springer-Verlag, 1973.

\bibitem {b-d-f-3} L.~Brown, R.~Douglas, and P.~Fillmore, {\it
Extensions of $C^*$-algebras and $K$-homology}, Annals of Mathematics
(2) {\bf 105} (1977), 265--324.

\bibitem {c-w-1} R.~Coifman and G.~Weiss, {\it Analyse Harmonique
Non-Commutative sur Certains Espaces Homog\`enes}, Lecture Notes in
Mathematics {\bf 242}, Springer-Verlag, 1971.

\bibitem {c-w-2} R.~Coifman and G.~Weiss, {\it Extensions of Hardy
spaces and their use in analysis}, Bulletin of the American
Mathematical Society {\bf 83} (1977), 569--645.

\bibitem {cns} A.~Connes, {\it Noncommutative Geometry}, Academic
Press, 1994.

\bibitem {d1} R.~Douglas, {\it $C^*$-Algebra Extensions and
$K$-Homology}, Princeton University Press, 1980.

\bibitem {d2} R.~Douglas, {\it Banach Algebra Techniques in Operator
Theory}, 2nd edition, Springer-Verlag, 1998.

\bibitem {du} P.~Duren, {\it Theory of $H^p$ Spaces}, Academic Press,
1970.

\bibitem {d-s} P.~Duren and A.~Schuster, {\it Bergman Spaces},
American Mathematical Society, 2004.

\bibitem {g1} J.~Garnett, {\it Analytic Capacity and Measure}, Lecture
Notes in Mathematics {\bf 297}, Springer-Verlag, 1972.

\bibitem {g2} J.~Garnett, {\it Bounded Analytic Functions},
Springer-Verlag, 2007.

\bibitem {k} J.~Kigami, {\it Analysis on Fractals}, Cambridge
University Press, 2001.

\bibitem {sk} S.~Krantz, {\it A Panorama of Harmonic Analysis},
Mathematical Association of America, 1999.

\bibitem {k-p} S.~Krantz and H.~Parks, {\it The Geometry of Domains in
Space}, Birkh\"auser, 1999.

\bibitem {p} V.~Peller, {\it Hankel Operators and their Applications},
Springer-Verlag, 2003.

\bibitem {sa} D.~Sarason, {\it Function Theory on the Unit Circle},
Virginia Polytechnic and State University, 1978.

\bibitem {st1} E.~Stein, {\it Singular Integrals and Differentiability
Properties of Functions}, Princeton University Press, 1970.

\bibitem {st2} E.~Stein, {\it Harmonic Analysis: Real-Variable
Methods, Orthogonality, and Oscillatory Integrals}, with the
assistance of T.~Murphy, Princeton University Press, 1993.

\bibitem {s-w} E.~Stein and G.~Weiss, {\it Introduction to Fourier
Analysis on Euclidean Spaces}, Princeton University Press, 1971.

\bibitem {str} R.~Strichartz, {\it Differential Equations on
Fractals}, Princeton University Press, 2006.




\end{thebibliography}
\end{document}